\documentclass[reqno, 11pt]{article}

\usepackage{enumerate}
\usepackage{needspace}
\usepackage{latexsym}
\usepackage{amsmath}
\usepackage{amsfonts}
\usepackage{amsthm}
\usepackage{amssymb,mathtools}
\usepackage{newlfont}
\usepackage{graphics}
\usepackage{color}
\usepackage{float}
\usepackage{diagbox}
\usepackage{tocloft}
\usepackage{titlesec}
\usepackage{booktabs,array,longtable,enumitem,microtype}
\usepackage{extpfeil}
\usepackage{centernot}

\usepackage[linesnumbered,ruled,vlined]{algorithm2e}
\usepackage{url}

\usepackage[margin=1in]{geometry}

\usepackage{tabularx}
\usepackage[T1]{fontenc}
\usepackage{lmodern}
\usepackage{rotating}
\usepackage{multirow}
\usepackage{extarrows}
\usepackage[sort,compress,numbers]{natbib}
\usepackage[utf8]{inputenc}
\usepackage{xcolor}
\usepackage{listings}
\usepackage{aliascnt}
\usepackage[pagebackref,colorlinks=true,linkcolor=blue,citecolor=red,urlcolor=blue]{hyperref}
\usepackage[nameinlink,noabbrev,capitalize]{cleveref}
\numberwithin{equation}{section}

\setlist[enumerate]{label=\textup{(\roman*)},leftmargin=*}
\newtheorem{theorem}{Theorem}[section]
\newaliascnt{proposition}{theorem}

\aliascntresetthe{proposition}
\newaliascnt{lemma}{theorem}
\newtheorem{lemma}[lemma]{Lemma}
\aliascntresetthe{lemma}
\newaliascnt{corollary}{theorem}
\newtheorem{corollary}[corollary]{Corollary}
\aliascntresetthe{corollary}

\theoremstyle{definition}
\newaliascnt{definition}{theorem}

\aliascntresetthe{definition}
\newaliascnt{remark}{theorem}

\aliascntresetthe{remark}
\newaliascnt{example}{theorem}

\aliascntresetthe{example}
\newaliascnt{conjecture}{theorem}

\aliascntresetthe{conjecture}
\newaliascnt{question}{theorem}

\aliascntresetthe{question}
\newaliascnt{problem}{theorem}

\aliascntresetthe{problem}
\newaliascnt{alg}{theorem}

\aliascntresetthe{alg}

\allowdisplaybreaks[4]

\SetKwInput{KwInput}{Input}
\SetKwInput{KwOutput}{Output}

\newcommand{\D}{\mathbb D}

\hypersetup{
  pdftitle={Proof of Stanley's Conjecture on the Leading Ehrhart Coefficients of the Birkhoff Polytope},
  pdfsubject={Ehrhart polynomials and enumeration of contingency tables},
  pdfdisplaydoctitle=true
}

\title{Proof of Stanley's Conjecture on the
Leading Ehrhart Coefficients of the Birkhoff Polytope}

\author{Jiahang Liu$^{\color{blue} \dag}$, Jinlong Tang$^{\color{blue} \ddag}$
\\[2mm]
{\small $^{\color{blue} \dag, \ddag}$ School of Mathematical Sciences,}\\[-0.8ex]
{\small Capital Normal University, Beijing, 100048, P.R.~China}\\
{\small {\color{blue} $^\dag$} Email address: \url{jiahang-liu@foxmail.com}}\\
{\small {\color{blue} $^\ddag$} Email address: \url{jinlong\_tang@cnu.edu.cn}}\\
}

\date{\today}

\begin{document}
\maketitle

\begin{abstract}
Let $c(n,k)$ denote the coefficient of $r^{(n-1)^2-k}$ in the
Ehrhart polynomial of the $n\times n$ Birkhoff polytope.
For every fixed nonnegative integer $k$, we prove that
$c(n,k)/c(n,0)\sim n^{3k}/(2^k k!)$ as $n\to\infty$, confirming
a conjecture of Stanley.
\end{abstract}

\noindent
\begin{small}
\end{small}

\noindent
\begin{small}
\end{small}

\section{Introduction}\label{sec:introduction}

For an integer $n\ge2$, the \emph{Birkhoff polytope} $B_n$ is the set of all nonnegative real $n\times n$ matrices whose row and column sums are all equal to $1$. For each nonnegative integer $r$, the \emph{Ehrhart function} $H_n(r)$ of $B_n$ counts the lattice points in the dilated polytope $rB_n$, namely,
\[
H_n(r)=\#\left(rB_n\cap\mathbb Z^{n^2}\right).
\]
Equivalently, $H_n(r)$ counts the nonnegative integer $n\times n$ matrices whose row and column sums are all equal to $r$. Such matrices are called \emph{semimagic squares}.

The Birkhoff polytope $B_n$ is an integral polytope of dimension $d=(n-1)^2$. Hence, by Ehrhart's theorem, its Ehrhart function is a polynomial of degree $d$:
\begin{equation}\label{eq:ehrhart}
H_n(r)=\sum_{k=0}^{d}c(n,k)r^{d-k}.
\end{equation}

The Ehrhart polynomial $H_n(r)$ connects the enumeration of semimagic squares with the geometry of $B_n$. Beck and Pixton~\cite{BP} developed formulas and computational methods for determining the Ehrhart polynomial and volume of $B_n$. In particular, their approach yields $H_9$ and the volume of $B_{10}$. More recently, Jiang, Xin, Zhang, and Zhong~\cite{JXZZ} computed the full Ehrhart series $\sum_{r\ge0}H_n(r)z^r$
for $n=10,\ldots,13$, thereby determining the corresponding Ehrhart polynomials using constant-term methods and finite fields. De Loera, Liu, and Yoshida~\cite{DLY} derived an explicit formula for each coefficient of $H_n$ as a sum indexed by permutations and rooted directed spanning trees on $n$ vertices. Canfield and McKay used their estimates for the number of integer matrices with equal margins~\cite{CM} to derive an asymptotic formula for the leading coefficient $c(n,0)$ as $n\to\infty$~\cite{CMV}.

We study the asymptotic behavior of $c(n,k)/c(n,0)$ for each fixed $k$ as $n\to\infty$. For the coefficient of the next-to-leading term, Stanley~\cite[Exercise~4.54(a) and its solution]{EC} records the exact formula
\[
\frac{c(n,1)}{c(n,0)}=\frac{n(n-1)^2}{2}.
\]
Stanley observed this identity empirically and credits its proof to
B.~Osserman and F.~Liu. In particular, the ratio for $k=1$ is asymptotic to $n^3/2$. In part~(b) of the same exercise, Stanley conjectured an asymptotic formula for $c(n,k)/c(n,0)$ for every fixed $k$. We establish this conjecture in the following form.

\begin{theorem}\label{thm:main}
For every fixed nonnegative integer $k$, we have
\[
 \frac{c(n,k)}{c(n,0)}
 \sim \frac{n^{3k}}{2^k k!}
 \qquad(n\to\infty).
\]
Here $c(n,k)$ is considered for all sufficiently large $n$, so that
$k\leq d$.
\end{theorem}

We introduce a sequence of auxiliary polynomials and establish their locally uniform convergence to an exponential function. This convergence yields the desired asymptotic formulas for the coefficient ratios. In \cref{sec:estimates}, we use Vitali's theorem to establish the locally uniform convergence of these polynomials on the open unit disk. We then use this convergence in \cref{sec:main-proof} to prove \cref{thm:main}.

\section{Convergence of the auxiliary polynomials}\label{sec:estimates}

Set
\begin{equation}\label{eq:normalization}
 N=n^3,\qquad A_n=\frac{e^{1/2}(n^2)!}{n^{d}\,(n!)^{2n}},
\end{equation}
and define the reversed, rescaled Ehrhart polynomial by
\begin{equation}\label{eq:rescaled-polynomial}
 F_n(z)=\sum_{k=0}^{d}\frac{c(n,k)}{A_nN^k}\,z^k.
\end{equation}
Each $F_n$ is holomorphic on the open unit disk $\D=\{z\in\mathbb C:|z|<1\}$. In this section, we prove that $F_n$ converges locally uniformly to $e^{z/2}$ on the open unit disk $\mathbb D$.

We use the following form of Vitali's theorem~\cite[Section~2.4]{Schiff}.

\begin{theorem}[Vitali convergence theorem]\label{Vi}
Let $D\subseteq\mathbb C$ be a domain, and let $(f_n)$ be a
sequence of holomorphic functions on $D$. Suppose that:
\begin{enumerate}
\item There exists a set $E\subset D$ with an accumulation point
in $D$ such that $(f_n(z))$ converges for every $z\in E$.
\item For every compact subset $K\subset D$, there exists a
constant $M_K>0$ such that
\[
|f_n(z)|\le M_K
\qquad\text{for all }z\in K\text{ and }n\ge1.
\]
\end{enumerate}
Then $(f_n)$ converges uniformly on every compact subset of $D$
to a holomorphic function $f$. Moreover, $f$ is uniquely determined
by its values on $E$.
\end{theorem}

\subsection{Pointwise convergence}\label{subsec:pointwise}

For $z\ne0$, reversing the powers in \eqref{eq:ehrhart} gives
\begin{equation}\label{eq:reversal}
 F_n(z)=\frac{z^d}{A_nN^d}\,H_n\!\left(\frac{N}{z}\right).
\end{equation}

At $z=1/q$ the count in \eqref{eq:reversal} is $H_n(qN)$.
We use the square case of Canfield and McKay's enumeration
estimate~\cite[Corollary~1]{CM}, with the hypotheses inherited
from their Theorem~1 made explicit.

\begin{lemma}[Canfield--McKay, square case]\label{CM}
Let $r=r(n)$ be a sequence of positive integers, and set
$\lambda=r/n$ and $d=(n-1)^2$.
Fix constants $a,b>0$ with $a+b<1/2$.
Suppose that
\begin{equation}\label{logn}
 \frac83\left(1+\frac{1}{4\lambda(1+\lambda)}\right)
 \leq a\log n
\end{equation}
for all sufficiently large $n$, and that
$n^2/\lambda^2\to0$ as $n\to\infty$.
Then
\begin{equation}\label{coro}
 H_n(r)
 =
 \left(\lambda+\frac12\right)^d
 \frac{(n^2)!}{(n!)^{2n}}
 \exp\left(
   \frac12+
   O\left(\frac{n^2}{\lambda^2}+n^{-b}\right)
 \right).
\end{equation}
\end{lemma}

This is obtained by setting $m=n$ and $s=t=r$ in
Canfield and McKay's Corollary~1.

\begin{lemma}\label{exp}
For every fixed positive integer $q$,
\[
 F_n(1/q)\longrightarrow \exp\left(\frac1{2q}\right)
 \qquad(n\to\infty).
\]
\end{lemma}

\begin{proof}
Fix a positive integer $q$ and take $r=qN=qn^3$ in \cref{CM}, so
that $\lambda=r/n=qn^2$. Choose $a=1/8$ and $b=1/4$, for which
$a+b<1/2$. Since $\lambda\to\infty$, the left-hand side of
\eqref{logn} tends to $8/3$, while $a\log n\to\infty$.
Thus \eqref{logn} holds for all sufficiently large $n$. The remaining hypothesis of \cref{CM} is also satisfied, since
\[
\frac{n^2}{\lambda^2}=\frac{1}{q^2n^2}\longrightarrow0.
\]

Using $qN=n\lambda$ and the definition of $A_n$ in \eqref{eq:normalization}, we divide both sides of \eqref{coro} by $A_n(qN)^d$. The identity \eqref{eq:reversal} then yields
\[
F_n(1/q)=\frac{H_n(qN)}{A_n(qN)^d}
=\left(1+\frac{1}{2qn^2}\right)^d
\exp\left(O\left(q^{-2}n^{-2}+n^{-1/4}\right)\right).
\]
For fixed $q$, we have
\[
\begin{aligned}
d\log\left(1+\frac{1}{2qn^2}\right)=\frac{d}{2qn^2}+O(q^{-2}n^{-2})\longrightarrow\frac{1}{2q},
\end{aligned}
\]
since $d/n^2\to1$. Thus the first factor converges to $e^{1/(2q)}$. Meanwhile, the exponential error factor tends to $1$. We combine these limits and get the stated result.
\end{proof}

\subsection{Local boundedness on the unit disk}\label{subsec:boundedness}

We show that the sequence $(F_n)_{n\geq2}$ is uniformly bounded on the closed unit disk.

\begin{lemma}\label{lem:product-comparison}
Let $a_1,\ldots,a_m\in\mathbb R$ satisfy
$$
\max_{1\le i\le m}|a_i|\le\frac12,
\qquad
\sum_{i=1}^m a_i^2\le1.
$$
Define $P(z)=\prod_{i=1}^m(1+a_i z)$. Then $P(1),P(-1)>0$. Furthermore, $P$ satisfies
$$
|P(z)|\le e^2\bigl(P(1)+P(-1)\bigr)
\qquad(|z|\le1).
$$
\end{lemma}

\begin{proof}
For $|w|\le1/2$, the power series of the principal logarithm gives
\[
 |\log(1+w)-w|
 \leq\sum_{\ell\geq2}\frac{|w|^\ell}{\ell}
 \leq\frac{|w|^2}{2(1-|w|)}
 \leq |w|^2.
\]
Every factor $1+a_i z$ lies in the open right half-plane for
$|z|\leq1$, so this logarithm is well-defined. Set $S=\sum_i a_i$
and $Q=\sum_i a_i^2\leq1$. Choose $\varepsilon\in\{-1,1\}$ so that
$\varepsilon S=|S|$. Taking real parts and summing gives
\[
 \log|P(z)|\leq(\operatorname{Re}z)S+|z|^2Q\leq |S|+Q,
 \qquad
 \log P(\varepsilon)\geq\varepsilon S-Q=|S|-Q.
\]
At the real endpoints all factors $1\pm a_i$ are positive. Thus
\[
 |P(z)|\leq e^{2Q}P(\varepsilon)
 \leq e^2\bigl(P(1)+P(-1)\bigr).
 \qedhere
\]
\end{proof}

\begin{lemma}\label{con}
The sequence $(F_n)_{n\geq2}$ is uniformly bounded on the closed unit
disk. That is, there exists a constant $C$, independent of $n$, such
that
\[
 \sup_{|z|\leq1}|F_n(z)|\leq C
\]
holds for all $n$.
\end{lemma}

\begin{proof}
Let $h_{n,0}^*,\ldots,h_{n,d}^*$ denote the $h^*$-coefficients of
$B_n$. Its Ehrhart polynomial has the binomial-basis expansion
\begin{equation}\label{eq:binomial-basis}
 H_n(r)=\sum_{j=0}^{d}h_{n,j}^*\binom{r+d-j}{d},
 \qquad
 \binom{r+d-j}{d}=\frac1{d!}\prod_{i=1}^{d}(r+i-j).
\end{equation}
Stanley's nonnegativity theorem~\cite[Theorem~2.1]{S} gives
$h_{n,j}^*\geq0$. The binomial coefficients appearing in \eqref{eq:binomial-basis} are polynomial functions of $r$. Hence both sides of the identity are polynomials in $r$, and their agreement for all nonnegative integers implies that the identity holds for every $r\in\mathbb C$.

Set $a_{ij}=(i-j)/N$ and define
\[
P_{n,j}(z)=\prod_{i=1}^{d}(1+a_{ij}z).
\]
By \eqref{eq:binomial-basis} and \eqref{eq:reversal}, we obtain
\begin{equation}\label{eq:positive-mixture}
F_n(z)=\sum_{j=0}^{d}w_{n,j}P_{n,j}(z),
\qquad
w_{n,j}=\frac{h_{n,j}^*}{d!A_n}\geq0.
\end{equation}
Since both sides are polynomials and agree for all $z\ne0$, the identity also holds at $z=0$. For $1\leq i\leq d$ and $0\leq j\leq d$, we have
$|i-j|\leq d$. Consequently, for every $n\geq2$, we get
\[
 \max_{1\leq i\leq d}|a_{ij}|
 \leq\frac dN\leq\frac1n\leq\frac12,\qquad
 \sum_{i=1}^{d}a_{ij}^2\leq\frac{d^3}{N^2}
 =\left(1-\frac1n\right)^6\leq1 .
\]
Applying \cref{lem:product-comparison} to each $P_{n,j}$ and summing
with the nonnegative weights in \eqref{eq:positive-mixture} gives
\begin{equation}\label{leq}
 \sup_{|z|\leq1}|F_n(z)|\leq e^2\bigl(F_n(1)+F_n(-1)\bigr).
 \end{equation}
By Ehrhart--Macdonald reciprocity~\cite[Theorem~4.6.9]{EC}, we have
\[
(-1)^d H_n(-N) = H_n^\circ(N),
\]
where $H_n^\circ(N)$ counts lattice points in the relative interior of $N B_n$. Since $H_n^\circ(N)\leq H_n(N)$, the reversal formula (\ref{eq:reversal}) yields
$$
0\leq F_n(-1)=\frac{(-1)^dH_n(-N)}{A_nN^d}=\frac{H_n^\circ(N)}{A_nN^d}
\leq\frac{H_n(N)}{A_nN^d}=F_n(1).
$$
It follows from \eqref{leq} that
\[
\sup_{|z|\le 1} |F_n(z)| \le 2e^2 F_n(1).
\]
Taking $q=1$ in \cref{exp} gives $F_n(1)\to e^{1/2}$, which
proves the claimed uniform bound.
\end{proof}

Now we set $E:=\{1/q:q=2,3,\ldots\}$. Then $E\subset \D$, and $0\in \D$ is an accumulation point of $E$. By \cref{exp}, $F_n(z)\to e^{z/2}$ for every $z\in E$, and
\cref{con} gives a uniform bound on $\D$. Since the $F_n$ are
holomorphic on $\D$ and $E$ has an accumulation point at $0\in\D$,
\cref{Vi} gives uniform convergence on compact subsets of $\D$ to a
holomorphic function $F$. Since $F$ agrees with $e^{z/2}$ on $E$,
the identity theorem identifies the limit throughout $\D$. Thus
\begin{equation}\label{eq:local-uniform-limit}
 F_n(z)\longrightarrow e^{z/2}
 \qquad\text{locally uniformly on $\D$}.
\end{equation}

\section{Proof of the main theorem}\label{sec:main-proof}

\begin{proof}[Proof of \cref{thm:main}]
Fix a nonnegative integer $k$ and choose $0<\rho<1$.
For all sufficiently large $n$, we have $k\leq d$. The convergence in
\eqref{eq:local-uniform-limit} is uniform on the circle $|z|=\rho$. Since the $k$-th Taylor coefficient of $e^{z/2}$ is
\[
[z^k] e^{z/2} = \frac{1}{2^k k!},
\]
Cauchy's coefficient formula then yields
\begin{align*}
 \left|[z^k]F_n(z)-\frac1{2^k k!}\right|
 &=\left|\frac1{2\pi i}\int_{|z|=\rho}
       \frac{F_n(z)-e^{z/2}}{z^{k+1}}\,dz\right|\\
 &\leq \rho^{-k}\sup_{|z|=\rho}|F_n(z)-e^{z/2}|
 \longrightarrow0.
\end{align*}
By \eqref{eq:rescaled-polynomial}, this yields
\begin{equation}\label{eq:coefficient-limit}
\frac{c(n,k)}{A_nN^k}=[z^k]F_n(z) \longrightarrow\frac1{2^k k!}.
\end{equation}
Evaluating \eqref{eq:local-uniform-limit} at $z=0$ gives
\begin{equation}\label{eq:volume-limit}
\frac{c(n,0)}{A_n}=F_n(0)\longrightarrow1.
\end{equation}
Recall that $N=n^3$. Dividing \eqref{eq:coefficient-limit} by
\eqref{eq:volume-limit}, we obtain
\[
\frac{c(n,k)}{c(n,0)n^{3k}}
=\frac{[z^k]F_n(z)}{F_n(0)}
\longrightarrow\frac1{2^k k!}.
\]
Hence $c(n,k)/c(n,0)\sim n^{3k}/(2^k k!)$, as claimed.
\end{proof}

The proof also yields the absolute asymptotics of the coefficients.

\begin{corollary}\label{cor:absolute-asymptotics}
For every fixed nonnegative integer $k$,
\[
 c(n,k)\sim
 \frac{e^{1/2}}{2^k k!}\,
 \frac{(n^2)!\,n^{3k}}
      {n^{(n-1)^2}(n!)^{2n}}.
\]
In particular, for each fixed $k$, the coefficient $c(n,k)$ is
positive for all sufficiently large $n$.
\end{corollary}

\begin{proof}
Substitute the definition of $A_n$ from \eqref{eq:normalization}
into \eqref{eq:coefficient-limit}.
\end{proof}

\vspace{5mm}
\noindent
\textbf{Acknowledgments.}
 The authors would like to thank the anonymous referees for valuable suggestions for improving the presentation.

\noindent\textbf{Declaration of AI assistance.}
During the preparation of this manuscript, the authors used ChatGPT to assist in exploring possible approaches, checking technical details, and improving the exposition. All mathematical arguments, computations, proofs, and results were independently verified by the authors. The authors take full responsibility for the content of this manuscript.

\end{document}